%% file: root.tex
\documentclass[letterpaper, 10 pt, conference]{ieeeconf}  % Comment this line out if you need a4paper

\IEEEoverridecommandlockouts                              % This command is only needed if 
\title{\LARGE \bf
Preparation of Papers for IEEE Sponsored Conferences \& Symposia*
}

\usepackage{amsmath,amssymb}

\usepackage{authblk}

\usepackage{color, colortbl}
\definecolor{Gray}{gray}{0.8}
\definecolor{LightCyan}{rgb}{0.88,1,1}
\definecolor{green(html/cssgreen)}{rgb}{0.0, 0.5, 0.0}

\newtheorem{remark}{Remark}
 
\newtheorem{theorem*}{Theorem}
\newtheorem{corollary}{Corollary}

\usepackage{mathtools}
\usepackage{mdframed}

\newmdtheoremenv{theo}{Theorem}

\input{preamble}

\title{Matrix Measure Flows: A Novel Approach \\ to Stable Plasticity in Neural Networks}

\author[1]{Leo Kozachkov}
\author[1,2]{Jean-Jacques Slotine}
\affil[1]{\small Department of Brain and Cognitive Sciences, Massachusetts Institute of Technology}
\affil[2]{\small Department of Mechanical Engineering, Massachusetts Institute of Technology}
\affil[ ]{\texttt{\{leokoz8,jjs\}@mit.edu}}
\date{}

\begin{document}

\maketitle
\thispagestyle{empty}
\pagestyle{empty}

%%%%%%%%%%%%%%%%%%%%%%%%%%%%%%%%%%%%%%%%%%%%%%%%%%%%%%%%%%%%%%%%%%%%%%%%%%%%%%%%
\begin{abstract}
This letter introduces the notion of a \textit{matrix measure flow} as a tool for analyzing the stability of neural networks with time-varying weights. Given a matrix flow -- for example, one induced by gradient-based adaptation -- the matrix measure flow tracks the progression of an associated matrix measure (or logarithmic norm). We show that for certain matrix flows of interest in computational neuroscience and machine learning, the associated matrix measure flow obeys a simple inequality. In the context of neuroscience, synapses -- the connections between neurons in the brain -- are constantly being updated. This plasticity subserves many important functions, such as memory consolidation and fast parameter adaptation towards real-time control objectives. However, synaptic plasticity poses a challenge for stability analyses of recurrent neural networks, which typically assume fixed synapses. Matrix measure flows allow the stability and contraction properties of recurrent neural networks with plastic synapses to be systematically analyzed. This in turn can be used to establish the robustness properties of neural dynamics, including those associated with problems in optimization and control. We consider examples of synapses subject to Hebbian and/or Anti-Hebbian plasticity, as well as covariance-based and gradient-based rules. 
\end{abstract}

\section{INTRODUCTION}
Biological synapses are always in flux. Changes to synaptic strength can occur on timescales between hundreds of milliseconds to hours and days \cite{dayan2005theoretical}. These changes underlie myraid important functions, such as long-term learning and fast biological adaptive control. As pointed out by \cite{zenke2017temporal}, this poses a challenge for the brain. Synaptic plasticity inside a recurrent neural network can lead to uncontrolled instabilities. This can readily be seen for correlation-based learning rules such as Hebbian learning. If two neurons have correlated activity, the synaptic weight between them increases. Since the weight increased, their activity will correlate more, thus increasing the weight again, and so on. This is a positive feedback loop that will lead to exponential ``explosions" of activity if left unchecked \cite{dayan2005theoretical}. Since this explosion does not appear to happen in the brain, it is of interest to ask ``what synaptic plasticity rules preserve overall stability"? Of course, this stability problem is not limited to biological brains. It applies broadly to machines with changing parameters.

\paragraph{Contraction Analysis} There has been much recent work in the computational neuroscience, machine learning, and control theory communities analyzing the stability of recurrent neural network with \textit{fixed} weights \cite{revay2020contracting,kozachkov2021rnns,davydov2022non,jafarpour2021robust,manchester2021contraction}. In particular, authors have been interested in restricting the weights in a recurrent network so that it is provably \textit{contracting} \cite{lohmiller1998contraction}. Any pair of trajectories from a contracting dynamical system will converge towards one another exponentially, regardless of initial condition. Contraction implies many desirable properties for neural networks and Neural ODEs, such as exponential convergence to a unique fixed point (for autonomous systems) \cite{slotine2003modular}, robustness \cite{jafarpour2021robust}, and the ability to be combined with other contracting system in ways that preserve contraction \cite{kozachkov2021rnns}. In this letter we extend these analyses to models with changing parameters, which is the more biologically realistic case. 

\paragraph{Neural Network} We consider an $n$-dimensional neural network (or Neural-ODE). The interactions between neurons in this network are governed by a single weight matrix $W$, through the following standard continuous-time equations
\begin{equation}\label{eq:RNN_eqs}
\tau\dot{x} = -x + W\phi(x) + u(t)    
\end{equation}
where $\tau \geq 0$ is the timescale of the neural network, $\phi$ is the nonlinearity, and $u(t)$ is some external input. Most applications of contraction analysis to \eqref{eq:RNN_eqs} have focused on the case where the weight matrix $W$ does not change. In this case, one can find simple restrictions on $W$ that guarantee global contraction of \eqref{eq:RNN_eqs}. However, an organism's ability to adapt to its environment crucial depends on $W$ changing. The ways in which the weight matrix changes is typically a function of the neural states, as well as the current weight matrix. We focus on weight updates which can be written in the form
\begin{equation}\label{eq:W_eqs}
\dot{W} = -\gamma W + G(W,x,t) 
\end{equation}
where $\gamma\geq 0$ is a leak coeffcient and $G$ is some nonlinear matrix function. Biologically, the leak term $-\gamma W$ serves the purpose of ensuring that synaptic weights do not ``blow up" \cite{gerstner2002mathematical}. It also emerges naturally in diverse settings, such as the minimization of Tikhonov regularized objective functions, as well as natural gradient descent in deep linear neural networks \cite{bernacchia2018exact}. The term $G$ captures interactions between the synaptic weights, neurons, as well as other external variables (e.g glia \cite{araque2014gliotransmitters}). For example, a simple correlational Hebbian rule  \cite{dong1992dynamic,dayan2005theoretical} can be written as the outer product of firing rate vector with itself
\[ G_{Hebb} = \phi(x)\phi(x)^T\]
Taken together, the equations \eqref{eq:RNN_eqs} and \eqref{eq:W_eqs} represent a closed-loop dynamical system consisting of two interacting components: a neural component and a synaptic components. In this note, we use contraction analysis to derive sufficient conditions on the learning rule $G$ which ensure that the combined systems \eqref{eq:RNN_eqs} and \eqref{eq:W_eqs} will remain stable.

\paragraph{Timescale Separation} We are interested in the case where the time constant of the neural network is small $\tau \ll 1$, or equivalently, when the dynamics of the synaptic weights are slower than those of the neurons. In this case, contraction of the overall system can be understood in the sense of singular perturbation theory, as explained in \cite{del2012contraction} and \cite{nguyen2018contraction}. Intuitively, because we assume a separation of timescales, from the perspective of the neural dynamics \eqref{eq:RNN_eqs} the synaptic weight matrix is constant. Conversely, from the perspective of the synaptic dynamics \eqref{eq:W_eqs}, the neural dynamics converge instantaneously. To maintain this interpretation, we need to ensure that the neural dynamics are contracting for all fixed $W$. As we will discuss in the next section, this can be achieved by ensuring the synaptic dynamics $\eqref{eq:W_eqs}$ force the \textit{matrix measure} of $W(t)$ to be sufficiently small as $t \rightarrow \infty$. 

\section{MATRIX MEASURES AND THEIR FLOWS}
\subsection{Matrix Measure Background}
The main tool we will use in our analysis is the matrix measure \cite{dahlquist1958stability,desoer1972measure}. Matrix measures have found wide application in dynamical systems and control theory, in particular in the context of contraction analysis. A matrix measure (or logarithmic norm) maps square matrices onto $\fR$. It is a function that is induced by a corresponding matrix norm $|| \cdot ||_i$, in the following way
\begin{equation}\label{eq:mm_def}
\mu_i\big[W\big] \coloneqq \lim_{h\to 0+} \frac{||I + hW||_i - 1}{h}
\end{equation}
Two common matrix measures are $\mu_2$ and $\mu_1$, the matrix measures associated with the matrix 2-norm and 1-norm, respectively
\begin{align*}
\mu_2\big[W\big] = \lambda_{max}\bigg[\frac{W + W^T}{2}\bigg] \\
\mu_1\big[W\big] = \max_j \bigg[W_{jj} + \sum_{i \neq j}^n |W_{ij}| \bigg] 
\end{align*}
where $\lambda_{max}[\cdot]$ denotes the largest eigenvalue. More details can be found in \cite{desoer1972measure} or \cite[Section 2.2.2]{vidyasagar2002nonlinear} . Two important properties of matrix measures $\mu \big[\cdot \big]$ are subadditivity and positive homogeneity. That is, for any two real matrices $A,B \in \fR^{n \times n}$, we have
\begin{align*}
\mu\big[A + B\big] \leq \mu\big[A\big] + \mu\big[B\big], &  &\text{(subadditivity)}  \\
\mu\big[cA\big]  = c\mu\big[A\big] \text{  for all  } c \geq 0, & &\text{(positive homogeneity)} 
\end{align*}
The reason matrix measures are important in contraction analysis is the following: if there exists a matrix measure of a system's Jacobian which is globally negative, then the system is globally contracting. For recurrent neural networks with fixed parameters, it is usually enough to restrict the matrix measure to ensure global contraction. These restrictions are usually upper bounds on the matrix measure, of the form\cite{fang1996stability,davydov2022non,kozachkov2021rnns}
\begin{equation}
\mu\big[W\big] \leq k \hspace{0.2cm} \text{with} \hspace{0.2cm} k > 0 
\end{equation}
The number $k$ is typically related to the nonlinearity $\phi$ used in \eqref{eq:RNN_eqs}. A common choice is $k = g^{-1}$ where $0 < \phi' \leq g$. For example, using this choice of $k$ it was shown by \cite{kozachkov2021rnns} that if $W$ is symmetric and $\mu_2\big[W \big] < g^{-1}$, then the RNN \eqref{eq:RNN_eqs} is globally contracting. For non-symmetric $W$, it was shown in \cite{matsuoka1992stability} that $\mu_2\big[W \big] < g^{-1}$ implies global asymptotic stability of a unique fixed point, provided that the inputs into the RNN are constant. Similarly, it was shown by \cite{davydov2022non} that $\mu_1\big[W \big] < g^{-1}$ is a sufficient condition for global contractivity of \eqref{eq:RNN_eqs}. 
\paragraph{Eigenvalues, Matrix Measures, and Matrix Norms} The stability of linear systems is fully characterized by the eigenvalues of the system matrix \cite{luenberger1979introduction}. Because nonlinear systems are sometimes well-approximated by linear systems around fixed points, eigenvalues have played a large role in the stability analysis of RNNs \cite{jaeger2001echo,christodoulou2022eigenvalue}. However, these analyses can typically only provide local stability results. Matrix norms and matrix measures have proven to be more useful in providing global stability results for the neural network \eqref{eq:RNN_eqs}. This is somewhat unsurprising, in view of the ordering between eigenvalues, matrix measures, and matrix norms. Note that for a given matrix norm $|| \cdot ||_i$ and matrix measure $\mu_i$, one has the following chain of inequalities related to the largest real part of the eigenvalues
\[\lambda_{max}\big[ W \big] \leq \mu_i\big[W\big] \leq || W ||_i \]
For nonlinear systems, eigenvalues tend to give liberal stability bounds, but these are only local. Matrix norms tend to give global bounds, but are often too conservative. Matrix measures span the Goldilocks zone in between these two, while still yielding global results.

\subsection{Matrix Measure Flows}
As explained above, we are interested in the situation when $W$ evolves with time according to \eqref{eq:W_eqs}. In this setting, we ask: what happens to the \textit{matrix measure} of $W$ when it evolves according to \eqref{eq:W_eqs}? The dynamics \eqref{eq:W_eqs} induce an unknown\textit{ matrix measure flow}
\[ \frac{d^+\mu\big[W \big]}{dt}  = \\ ? \]
where $\frac{d^+}{dt}$ denotes the Dini derivative \cite{khalil2015nonlinear} (or upper right-hand derivative) at time $t$. The reason we consider the Dini derivative--as opposed to an ordinary derivative--is because the matrix measure may not always be differentiable (consider the $\mu_1$ matrix measure above). It is, however, always continuous with respect to the matrix elements, which is all we need. In the following Theorem, we show that this matrix measure flow is related in a simple way to the evolution of $W$ itself.

\begin{theo}\label{theorem:mm_subadditive}
Let $\mu \big[\cdot \big]: \fR^{n \times n} \rightarrow \fR$ denote a matrix measure. Assume that $W \in \fR^{n \times n}$ evolves according to the matrix differential equation \eqref{eq:W_eqs}. Then the following scalar differential inequality holds
\begin{equation}\label{eq:mu_diff_ineq}
\frac{d^+\mu\big[W\big]}{dt} \leq -\gamma \mu\big[W\big] + \mu\big[G(W,x,t)\big]
\end{equation}
\end{theo}

\begin{proof}
To avoid clutter, we will write $G(W_t,x_t,t)$ as $G_t$. Consider the first-order Euler discretization of the matrix differential equation \eqref{eq:W_eqs}, with step-size $h$ such that $ \ h < \frac{1}{\gamma} $ ,

\[W_{t+h} = W_t - h\gamma W_t + h G_t  = (1-h\gamma)W_t + h G_t \]
Using the sub-additivity property of matrix measures, we have
\[\mu\big[W_{t+h}\big] = \mu\big[(1-h\gamma)W_t + h G_t\big] \leq  (1-h\gamma)\mu\big[W_t\big] + h\mu\big[G_t\big]\]
where the last step uses the positive homogeneity of matrix measures to `pull out' the $ \ 1-h\gamma \ $ and $ h $ terms. Subtracting $\mu\big[W_t\big]$ from both sides of this inequality and dividing both sides by $h$, we see that
\[\frac{\mu\big[W_{t+h}\big] - \mu\big[W_t\big]}{h} \leq -\gamma \mu\big[W_t\big]  + \mu\big[G_t\big]\]
Since the matrix measure $\mu\big[W_{t}\big]$ may not always be differentiable with respect to time, we consider the right-hand supremum limit (i.e., the Dini derivative) to arrive at the stated differential inequality
\begin{equation*}
\resizebox{1\hsize}{!}{$\frac{d^+ \mu\big[W\big]}{dt} = \limsup\limits_{h \rightarrow 0+} \frac{\mu\big[W_{t+h}\big] - \mu\big[W_t\big]}{h}  \leq -\gamma \mu\big[W\big] + \mu\big[G(W,x,t)\big]$ }
\end{equation*}
\end{proof}
\begin{remark}
Theorem \ref{theorem:mm_subadditive} holds even when the leak rate $\gamma$ depends on $W$ and $t$, so long as $\gamma(W,t) \geq 0$. 
\end{remark}

\begin{remark}
Theorem \ref{theorem:mm_subadditive} is applicable for \textit{any} function $\mu \big[\cdot \big]: \fR^{n \times n} \rightarrow \fR$ such that $\mu \big[\cdot \big]$ is subadditive and positively homogeneous. This includes, for example, any matrix norm. 
\end{remark}

\begin{corollary}\label{corollary: boundedness}
Assume $\gamma > 0$ and $\mu \big[G(W,x,t) \big] \leq D$, for some $D \in \fR$. Then equation \eqref{eq:mu_diff_ineq} implies an upper-bound on $\mu \big[ W(t)\big]$

\begin{equation}\label{equation: upper-bound-mu}
\mu \big[ W(t) \big] \leq \mu \big[W(0)\big]e^{-\gamma t} + \frac{D}{\gamma}     
\end{equation}
This is a consequence of the Gr\"{o}nwall Comparison Lemma (see Exercise E2.1 in \cite{FB-CTDS} for an explicit proof). From \eqref{equation: upper-bound-mu} we can see that the inequality $\mu \big[ W(t) \big] \leq k$ will be satisfied in finite time if the following strict inequality between $D,k$ and $\gamma$ is satisfied
\begin{equation}\label{equation: strict_ineq_D_k}
\frac{D}{k} < \gamma
\end{equation}
This strict inequality is necessary to ensure that $\mu \big[ W(t) \big]$ eventually drops below $k$, instead of merely approaching it from above. 
\end{corollary}

\section{APPLICATIONS}
Theorem \ref{theorem:mm_subadditive} has immediate application to the problem of determining contraction in neural networks with synaptic plasticity.
\subsection{Anti-Hebbian Plasticity}\label{subsection: anti-hebb}
Consider the simple anti-Hebbian plasticity rule analyzed in \cite{kozachkov2020achieving} and \cite{centorrino2023euclidean}, with
\[\dot{W} = -\gamma W -\phi(x)\phi(x)^T \]
Using the notation of \eqref{eq:W_eqs}, we can identify
\[G(W,x,t) = -\phi(x)\phi(x)^T \]
Since $G(W,t)$ is a negative semi-definite matrix, we have that $\mu_2\big[G\big] \leq 0$. This implies that $D$ in Corollary \ref{corollary: boundedness} may be taken to be $D = 0$. Substituting $D = 0$ into equation \eqref{equation: upper-bound-mu}, we see that $\mu_2\big[W(t)\big]$ decays exponentially to zero
\[ \mu_2\big[W(t)\big] \leq \mu_2\big[W(0)\big]e^{-\gamma t} \]
Note that the above exponential decay does \textit{not} imply that $W$ itself decays to zero. It simply implies that largest eigenvalue of the symmetric part of $W$ decays towards zero. Since $G$ is symmetric, we can also deduce that if $W(t)$ is initialized as a symmetric matrix, then $W(t)$ will be symmetric for all time. If $W$ is initially asymmetric, then it will converge exponentially with rate $\gamma$ towards a symmetric matrix. This follows from decomposing $W$ into a symmetric and skew-symmetric matrix, and then noting that the skew-symmetric part decays exponentially with rate $\gamma$.     
\begin{figure*}[ht]
\label{fig: fig1}
\centering
\includegraphics[width=\textwidth]{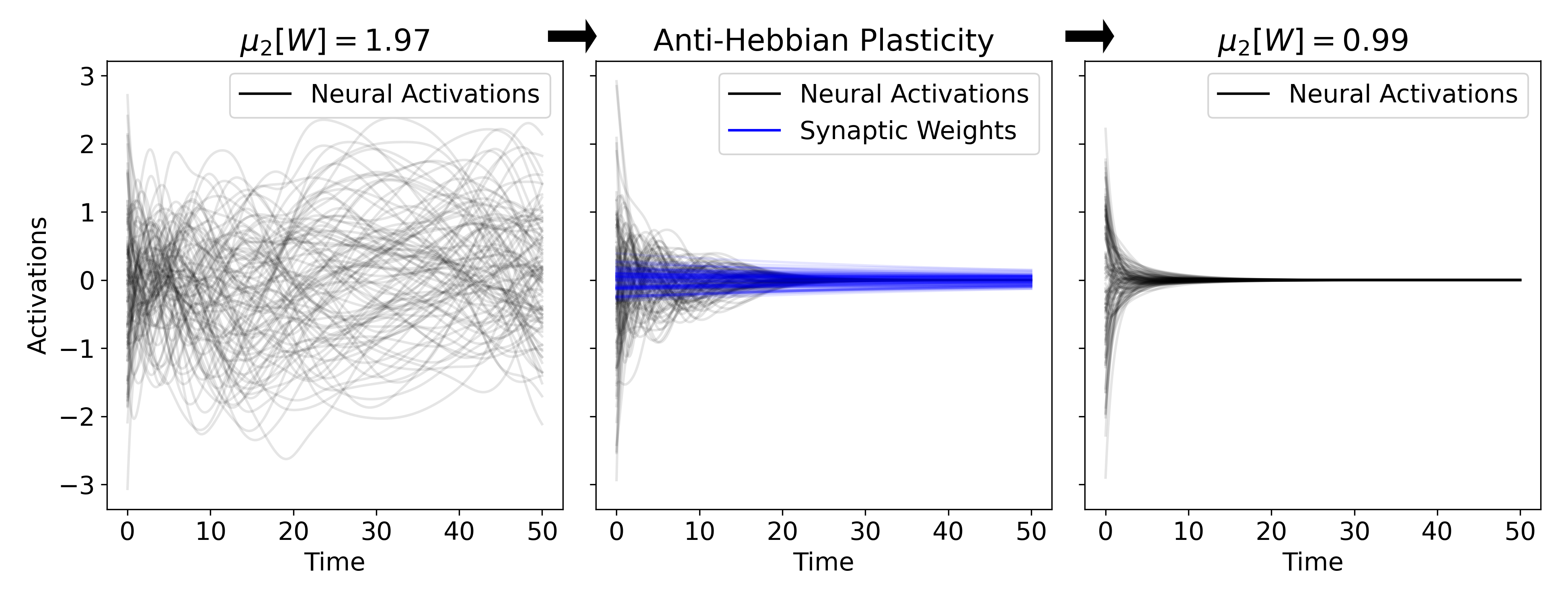}
\caption{Left) Chaotic neural dynamics for the RNN described by \eqref{eq:RNN_eqs}, associated with a large matrix measure for the synaptic weight matrix. Each trace the time-dependent activity of a neuron, as the network evolves from a randomly chosen initial condition. Center) Anti-Hebbian plasticity induces a stabilizing matrix measure flow. This flow forces the matrix measure of the weight matrix to be less than unity, leading to stable neural dynamics. See subsection \ref{subsection: hebb-anti-hebb} for details. Right) After the synaptic dynamics are run for a sufficiently long (yet finite) time, the weight matrix is such that the associated neural dynamics are stable. In this case the neural dynamics converge exponentially towards a unique, global attractor.}
\end{figure*}

\subsection{Hebbian/Anti-Hebbian Synaptic Plasticity}\label{subsection: hebb-anti-hebb}
The brain has both Hebbian and Anti-Hebbian synaptic plasticity \cite{bell1993storage,gerstner2002mathematical,pehlevan2015hebbian}. Theorem \ref{theorem:mm_subadditive} gives simple sufficient conditions for what combinations of Hebbian and Anti-Hebbian plasticity lead to contracting weight matrices. In particular, consider
\[\dot{W} = -K \odot \phi(x)\phi(x)^T -\gamma W \]
where $K \in \fR^{n \times n}$ is an arbitrary positive semi-definite matrix, and $\odot$ denotes element-wise matrix multiplication (i.e. the Hadamard product). These dynamics were also initially analyzed in \cite{kozachkov2020achieving}, but there the signs of the entries of $K$ were constrained to be positive. Here we do not constrain the signs of the individual matrix entries of $K$, so this plasticity rule is both Hebbian and/or Anti-Hebbian.  In this case, we have
\[G(W,t) = -K \odot \phi(x)\phi(x)^T \]
Since $-K$ is negative semi-definite by assumption and $\phi(x)\phi(x)^T$ is positive semi-definite by construction, the Schur Product Theorem \cite{zhang2006schur} tells us that their elementwise product will be negative semi-definite
\[\mu_2\big[-K \odot \phi(x)\phi(x)^T] \leq 0\]
This again implies that $D = 0$. As in the previous section, this in turn implies that the largest eigenvalue of the symmetric part of $W$ will decay to zero exponentially. Additionally, the symmetry of $K$ implies that $W$ will remain symmetric if initialized as a symmetric matrix, or will converge exponentially with rate $\gamma$ towards a symmetric matrix if initialized as an asymmetric matrix. 

Note that we may also include a matrix bias term $B = B^T$ in the synaptic dynamics
\[G(W,t) = -K \odot \phi(x)\phi(x)^T + B\]
From the sub-additivity of the matrix measure, the inclusion of $B$ increases the disturbance term $D$ by at worst $\mu\big[D\big]$. See Figure 1 for a numerical example, with $n = 100$, $\phi(x) = \text{tanh}(x)$, $K_{ij} = 1$ and $B_{ij} = 0.9/n$.

\subsection{Dong \& Hopfield Synaptic Plasticity}
The previous two examples relied on the positive semi-definiteness of the system matrix $K$, which in turn implies that $D \leq 0$. In this example we consider a case where $D \geq 0$. Consider the following Hebbian synaptic dynamics, first introduced and analyzed by \cite{dong1992dynamic}
\[\dot{W} = \nu\phi(x)\phi(x)^T -\gamma W \]
where $\nu > 0$ is a learning rate. Consider for example the bounded activation $\phi(\cdot) = \text{tanh}(x)$. Applying the matrix measure associated with the 2-norm to
\[G(W,t) = \nu\phi(x)\phi(x)^T \]
we can see that $D \leq \nu\phi(x)^T\phi(x) \leq \nu n$. This follows from the observation that the only non-zero eigenvalue of $G$ is $\nu\phi(x)^T\phi(x)$, and that $\phi$ is bounded between negative one and one. From this we can conclude that if
\[ \frac{\nu n}{k} < \gamma\] 
the weight matrix will satisfy the inequality $\mu_2 \big[W(t)\big] \leq k$ in finite time.

\subsection{Covariance Rule}
Consider a covariance based learning rule \cite{sejnowski1989hebb,gerstner2002mathematical}, with a leak term
\[ \dot{W} = \nu \big[ \phi(x) - \langle \phi \rangle  \big]\big[ \phi(x) - \langle \phi \rangle   \big] ^T - \gamma W \]
where $\nu \geq 0$ is a learning rate and $\langle \phi \rangle $ denotes a running average over the firing rate
\[\langle \phi \rangle_i = \int_{t-\Delta}^{t} \phi(x(s))_i \,ds\]
with some window size $\Delta > 0$. As in the section above, we can identity $D = \nu||\big[ \phi(x) - \langle \phi \rangle ||^2$ using the matrix measure associated with the 2-norm. If we assume that the firing rate does not deviate more than $\sigma^2$ from its running average, then the strict inequality
\[\frac{\nu\sigma^2}{k} < \gamma \]
implies that the weight matrix will satisfy the inequality $\mu_2 \big[W(t)\big] \leq k$ in finite time.

\subsection{Presynaptic Plasticity}
In some models \cite{mongillo2008synaptic,masse2019circuit}, synaptic plasticity only depends on the activity of the presynaptic neuron
\[\dot{W}_{ij} = -\gamma W_{ij} + b_i\phi(x_j) \]
where $b_i$ is some gain parameter. If we assume that the firing rate is bounded, such that $|\phi(x_j)| \leq \phi_{max}$, and the gain term is bounded such that $|b_i| \leq b_{max}$ then application of the matrix measure associated with the 1-norm yields the strict inequality
\[ \frac{nb_{max}\phi_{max}}{k} < \gamma   \]

\subsection{Gradient Descent} 
Many loss functions used in machine learning include a regularization term to control the model complexity. The total loss minimized during training is the sum of a ``task" loss term, which measures how poorly the model is doing, and a regularization term
\[\mathcal{L} = \mathcal{L}_{task} + \frac{\gamma}{2} \mathcal{L}_{reg} \]
where $\mathcal{L}_{task}$ denotes the loss purely due to performance error, and $\mathcal{L}_{reg}$ denotes the regularization term. A common choice of regularization is the squared sum-of-squares of all the weights (i.e the squared Frobenius norm). For square matrices this can be written as
\[\mathcal{L}_{reg} = \text{Tr}(W^T W) = ||W||^2_{fro}\]
Gradient flow on this loss function yields the following familiar expression for $W$
\begin{align*}
\dot{W} = -\nabla_W \mathcal{L} =  -\nabla_W \mathcal{L}_{task} - \frac{\gamma}{2} \nabla_W \mathcal{L}_{reg} \\ = -\nabla_W \mathcal{L}_{task} - \gamma W    
\end{align*}
If we assume that the loss function is Lipschitz with respect to the parameters $W$, with Lipschitz constant $L$, then it is straightforward to show that making $\gamma$ large enough will produce a weight matrix with a matrix measure less than $k$. Using the fact that the matrix measure is always upper-bounded by its associated matrix norm, we have
\[ \mu_i \big[ G \big] = \mu_i \big[ -\nabla_W \mathcal{L}_{task} \big] \leq ||-\nabla_W \mathcal{L}_{task}||_i \leq L \]
This shows that we may take $D$ in Corollary \eqref{corollary: boundedness} to be equal to $D = L$. To satisfy the inequality for $k$, we simply need to ensure that
\[\frac{L}{k}  < \gamma \]
\section{CONCLUSION}
Understanding how the brain balances stability and plasticity is an old and important problem in computational neuroscience \cite{zenke2017temporal}. We introduce the notion of a matrix measure flow and show that it is useful for approaching this problem. While our analysis has been focused on recurrent neural networks, the basic idea of a matrix measure flow can be applied to any system with ``plastic" weights (for example, biochemical reactions characterized by changing rate parameters \cite{russo2010global,raveh2016model}). Similarly, our results may be useful in the context of adaptive control and prediction, where certain adaptive laws may be thought of as gradient flow \cite{boffi2021random}. In this case, the assumption of timescale separation may be relaxed.

\addtolength{\textheight}{-12cm}   % This command serves to balance the column lengths

%%%%%%%%%%%%%%%%%%%%%%%%%%%%%%%%%%%%%%%%%%%%%%%%%%%%%%%%%%%%%%%%%%%%%%%%%%%%%%%%

%%%%%%%%%%%%%%%%%%%%%%%%%%%%%%%%%%%%%%%%%%%%%%%%%%%%%%%%%%%%%%%%%%%%%%%%%%%%%%%%

%%%%%%%%%%%%%%%%%%%%%%%%%%%%%%%%%%%%%%%%%%%%%%%%%%%%%%%%%%%%%%%%%%%%%%%%%%%%%%%%

\section*{ACKNOWLEDGMENTS}
We would like to thank Francesco Bullo for helpful suggestions and comments. L.K. is supported by a K. Lisa Yang Integrative and Computational Neuroscience (ICoN) Postdoctoral Fellowship.

\bibliographystyle{ieeetran}
\bibliography{mmf}

\end{document}

%% file: preamble.tex
\newcommand{\fR}{\mathbb{R}}

\usepackage[most]{tcolorbox}

%% file: root.bbl
\begin{thebibliography}{10}
\providecommand{\url}[1]{#1}
\csname url@rmstyle\endcsname
\providecommand{\newblock}{\relax}
\providecommand{\bibinfo}[2]{#2}
\providecommand\BIBentrySTDinterwordspacing{\spaceskip=0pt\relax}
\providecommand\BIBentryALTinterwordstretchfactor{4}
\providecommand\BIBentryALTinterwordspacing{\spaceskip=\fontdimen2\font plus
\BIBentryALTinterwordstretchfactor\fontdimen3\font minus
  \fontdimen4\font\relax}
\providecommand\BIBforeignlanguage[2]{{%
\expandafter\ifx\csname l@#1\endcsname\relax
\typeout{** WARNING: IEEEtran.bst: No hyphenation pattern has been}%
\typeout{** loaded for the language `#1'. Using the pattern for}%
\typeout{** the default language instead.}%
\else
\language=\csname l@#1\endcsname
\fi
#2}}

\bibitem{dayan2005theoretical}
P.~Dayan and L.~F. Abbott, \emph{Theoretical neuroscience: computational and
  mathematical modeling of neural systems}.\hskip 1em plus 0.5em minus
  0.4em\relax MIT press, 2005.

\bibitem{zenke2017temporal}
F.~Zenke, W.~Gerstner, and S.~Ganguli, ``The temporal paradox of hebbian
  learning and homeostatic plasticity,'' \emph{Current opinion in
  neurobiology}, vol.~43, pp. 166--176, 2017.

\bibitem{revay2020contracting}
M.~Revay and I.~Manchester, ``Contracting implicit recurrent neural networks:
  Stable models with improved trainability,'' in \emph{Learning for Dynamics
  and Control}.\hskip 1em plus 0.5em minus 0.4em\relax PMLR, 2020, pp.
  393--403.

\bibitem{kozachkov2021rnns}
L.~Kozachkov, M.~M. Ennis, and J.-J. Slotine, ``Rnns of rnns: Recursive
  construction of stable assemblies of recurrent neural networks,'' in
  \emph{Advances in Neural Information Processing Systems}, 2021.

\bibitem{davydov2022non}
A.~Davydov, A.~V. Proskurnikov, and F.~Bullo, ``Non-euclidean contractivity of
  recurrent neural networks,'' in \emph{2022 American Control Conference
  (ACC)}.\hskip 1em plus 0.5em minus 0.4em\relax IEEE, 2022, pp. 1527--1534.

\bibitem{jafarpour2021robust}
S.~Jafarpour, A.~Davydov, A.~Proskurnikov, and F.~Bullo, ``Robust implicit
  networks via non-euclidean contractions,'' \emph{Advances in Neural
  Information Processing Systems}, vol.~34, pp. 9857--9868, 2021.

\bibitem{manchester2021contraction}
I.~R. Manchester, M.~Revay, and R.~Wang, ``Contraction-based methods for stable
  identification and robust machine learning: a tutorial,'' in \emph{2021 60th
  IEEE Conference on Decision and Control (CDC)}.\hskip 1em plus 0.5em minus
  0.4em\relax IEEE, 2021, pp. 2955--2962.

\bibitem{lohmiller1998contraction}
W.~Lohmiller and J.-J.~E. Slotine, ``On contraction analysis for non-linear
  systems,'' \emph{Automatica}, vol.~34, no.~6, pp. 683--696, 1998.

\bibitem{slotine2003modular}
J.-J.~E. Slotine, ``Modular stability tools for distributed computation and
  control,'' \emph{International Journal of Adaptive Control and Signal
  Processing}, vol.~17, no.~6, pp. 397--416, 2003.

\bibitem{gerstner2002mathematical}
W.~Gerstner and W.~M. Kistler, ``Mathematical formulations of hebbian
  learning,'' \emph{Biological cybernetics}, vol.~87, no.~5, pp. 404--415,
  2002.

\bibitem{bernacchia2018exact}
A.~Bernacchia, M.~Lengyel, and G.~Hennequin, ``Exact natural gradient in deep
  linear networks and its application to the nonlinear case,'' \emph{Advances
  in Neural Information Processing Systems}, vol.~31, 2018.

\bibitem{araque2014gliotransmitters}
A.~Araque, G.~Carmignoto, P.~G. Haydon, S.~H. Oliet, R.~Robitaille, and
  A.~Volterra, ``Gliotransmitters travel in time and space,'' \emph{Neuron},
  vol.~81, no.~4, pp. 728--739, 2014.

\bibitem{dong1992dynamic}
D.~W. Dong and J.~J. Hopfield, ``Dynamic properties of neural networks with
  adapting synapses,'' \emph{Network: Computation in Neural Systems}, vol.~3,
  no.~3, p. 267, 1992.

\bibitem{del2012contraction}
D.~Del~Vecchio and J.-J.~E. Slotine, ``A contraction theory approach to
  singularly perturbed systems,'' \emph{IEEE Transactions on Automatic
  Control}, vol.~58, no.~3, pp. 752--757, 2012.

\bibitem{nguyen2018contraction}
H.~D. Nguyen, T.~L. Vu, K.~Turitsyn, and J.-J. Slotine, ``Contraction and
  robustness of continuous time primal-dual dynamics,'' \emph{IEEE control
  systems letters}, vol.~2, no.~4, pp. 755--760, 2018.

\bibitem{dahlquist1958stability}
G.~Dahlquist, ``Stability and error bounds in the numerical integration of
  ordinary differential equations,'' Ph.D. dissertation, Almqvist \& Wiksell,
  1958.

\bibitem{desoer1972measure}
C.~Desoer and H.~Haneda, ``The measure of a matrix as a tool to analyze
  computer algorithms for circuit analysis,'' \emph{IEEE Transactions on
  Circuit Theory}, vol.~19, no.~5, pp. 480--486, 1972.

\bibitem{vidyasagar2002nonlinear}
M.~Vidyasagar, \emph{Nonlinear systems analysis}.\hskip 1em plus 0.5em minus
  0.4em\relax SIAM, 2002.

\bibitem{fang1996stability}
Y.~Fang and T.~G. Kincaid, ``Stability analysis of dynamical neural networks,''
  \emph{IEEE Transactions on Neural Networks}, vol.~7, no.~4, pp. 996--1006,
  1996.

\bibitem{matsuoka1992stability}
K.~Matsuoka, ``Stability conditions for nonlinear continuous neural networks
  with asymmetric connection weights,'' \emph{Neural networks}, vol.~5, no.~3,
  pp. 495--500, 1992.

\bibitem{luenberger1979introduction}
D.~G. Luenberger, \emph{Introduction to dynamic systems: theory, models, and
  applications}.\hskip 1em plus 0.5em minus 0.4em\relax Wiley New York, 1979,
  vol.~1.

\bibitem{jaeger2001echo}
H.~Jaeger, ``The “echo state” approach to analysing and training recurrent
  neural networks-with an erratum note,'' \emph{Bonn, Germany: German National
  Research Center for Information Technology GMD Technical Report}, vol. 148,
  no.~34, p.~13, 2001.

\bibitem{christodoulou2022eigenvalue}
G.~Christodoulou and T.~Vogels, ``The eigenvalue value (in neuroscience),''
  \emph{OSF Preprints}, 2022.

\bibitem{khalil2015nonlinear}
H.~K. Khalil, \emph{Nonlinear control}.\hskip 1em plus 0.5em minus 0.4em\relax
  Pearson New York, 2015, vol. 406.

\bibitem{FB-CTDS}
\BIBentryALTinterwordspacing
F.~Bullo, \emph{Contraction Theory for Dynamical Systems}, {1.0}~ed.\hskip 1em
  plus 0.5em minus 0.4em\relax Kindle Direct Publishing, 2022. [Online].
  Available: \url{http://motion.me.ucsb.edu/book-ctds}
\BIBentrySTDinterwordspacing

\bibitem{kozachkov2020achieving}
L.~Kozachkov, M.~Lundqvist, J.-J. Slotine, and E.~K. Miller, ``Achieving stable
  dynamics in neural circuits,'' \emph{PLoS computational biology}, vol.~16,
  no.~8, p. e1007659, 2020.

\bibitem{centorrino2023euclidean}
V.~Centorrino, A.~Gokhale, A.~Davydov, G.~Russo, and F.~Bullo, ``Euclidean
  contractivity of neural networks with symmetric weights,'' \emph{IEEE Control
  Systems Letters}, 2023.

\bibitem{bell1993storage}
C.~C. Bell, A.~Caputi, K.~Grant, and J.~Serrier, ``Storage of a sensory pattern
  by anti-hebbian synaptic plasticity in an electric fish.'' \emph{Proceedings
  of the National Academy of Sciences}, vol.~90, no.~10, pp. 4650--4654, 1993.

\bibitem{pehlevan2015hebbian}
C.~Pehlevan, T.~Hu, and D.~B. Chklovskii, ``A hebbian/anti-hebbian neural
  network for linear subspace learning: A derivation from multidimensional
  scaling of streaming data,'' \emph{Neural computation}, vol.~27, no.~7, pp.
  1461--1495, 2015.

\bibitem{zhang2006schur}
F.~Zhang, \emph{The Schur complement and its applications}.\hskip 1em plus
  0.5em minus 0.4em\relax Springer Science \& Business Media, 2006, vol.~4.

\bibitem{sejnowski1989hebb}
T.~J. Sejnowski and G.~Tesauro, ``The hebb rule for synaptic plasticity:
  algorithms and implementations,'' in \emph{Neural models of
  plasticity}.\hskip 1em plus 0.5em minus 0.4em\relax Elsevier, 1989, pp.
  94--103.

\bibitem{mongillo2008synaptic}
G.~Mongillo, O.~Barak, and M.~Tsodyks, ``Synaptic theory of working memory,''
  \emph{Science}, vol. 319, no. 5869, pp. 1543--1546, 2008.

\bibitem{masse2019circuit}
N.~Y. Masse, G.~R. Yang, H.~F. Song, X.-J. Wang, and D.~J. Freedman, ``Circuit
  mechanisms for the maintenance and manipulation of information in working
  memory,'' \emph{Nature neuroscience}, vol.~22, no.~7, pp. 1159--1167, 2019.

\bibitem{russo2010global}
G.~Russo, M.~Di~Bernardo, and E.~D. Sontag, ``Global entrainment of
  transcriptional systems to periodic inputs,'' \emph{PLoS computational
  biology}, vol.~6, no.~4, p. e1000739, 2010.

\bibitem{raveh2016model}
A.~Raveh, M.~Margaliot, E.~D. Sontag, and T.~Tuller, ``A model for competition
  for ribosomes in the cell,'' \emph{Journal of The Royal Society Interface},
  vol.~13, no. 116, p. 20151062, 2016.

\bibitem{boffi2021random}
N.~M. Boffi, S.~Tu, and J.-J.~E. Slotine, ``Random features for adaptive
  nonlinear control and prediction,'' \emph{arXiv preprint arXiv:2106.03589},
  2021.

\end{thebibliography}
